\documentclass[12pt,a4paper,twoside]{article}

\newcommand{\pageformat}[6]{\setlength{\hoffset}{-1in}
                  \setlength{\voffset}{-1in}
                  \addtolength{\hoffset}{#5}
                            \addtolength{\voffset}{#6}
                            \setlength{\oddsidemargin}{#1}
                            \setlength{\evensidemargin}{#2}
                            \setlength{\textwidth}{\paperwidth}
                  \addtolength{\textwidth}{-\oddsidemargin}
                  \addtolength{\textwidth}{-\evensidemargin}
                  \addtolength{\textwidth}{-\marginparsep}
                  \addtolength{\textwidth}{-\marginparwidth}
                            \setlength{\topmargin}{#3}
                            \setlength{\textheight}{\paperheight}
                  \addtolength{\textheight}{-\topmargin}
                  \addtolength{\textheight}{-\headheight}
                  \addtolength{\textheight}{-\headsep}
                  \addtolength{\textheight}{-\footskip}
                  \addtolength{\textheight}{-#4}}
\pageformat{2cm}{3cm}{25mm}{25mm}{1pt}{0pt}

\usepackage{ifthen}
\newboolean{article}
    \setboolean{article}{true}
\newboolean{report}
\newboolean{book}
\newboolean{letter}
\newboolean{german}
\newboolean{nobaselinestretch}
\newboolean{nosectionappendix}
\newboolean{oldtoc}
\newboolean{nosectionequn}
\newboolean{notheorem}

\ifthenelse{\boolean{german}}{
    \usepackage{german}}{}

\usepackage[latin1]{inputenc}

\usepackage{amsmath}
\usepackage{amssymb}
\usepackage[mathscr]{eucal}

\ifthenelse{\boolean{notheorem}}{}{
    \usepackage{theorem}}

\ifthenelse{\boolean{nobaselinestretch}}{}{
    \renewcommand{\baselinestretch}{1.25}}

\newenvironment{env}[2]{\begin{#1}#2\end{#1}}{}
    \newcommand{\beq}[1]{\begin{env}{equation}{#1}}
    \newcommand{\beqn}[1]{\begin{env}{equation*}{#1}}
    \newcommand{\bal}[1]{\begin{env}{align}{#1}}
    \newcommand{\baln}[1]{\begin{env}{align*}{#1}}
    \newcommand{\bga}[1]{\begin{env}{gather}{#1}}
    \newcommand{\bgan}[1]{\begin{env}{gather*}{#1}}
    \newcommand{\bflal}[1]{\begin{env}{flalign}{#1}}
    \newcommand{\bflaln}[1]{\begin{env}{flalign*}{#1}}
    \newcommand{\bmu}[1]{\begin{env}{multline}{#1}}
    \newcommand{\bmun}[1]{\begin{env}{multline*}{#1}}
    \newcommand{\bsp}[1]{\begin{env}{split}{#1}}

    \newcommand{\eeq}{\end{env}}
    \newcommand{\eeqn}{\end{env}}
    \newcommand{\eal}{\end{env}}
    \newcommand{\ealn}{\end{env}}
    \newcommand{\ega}{\end{env}}
    \newcommand{\egan}{\end{env}}
    \newcommand{\eflal}{\end{env}}
    \newcommand{\eflaln}{\end{env}}
    \newcommand{\emu}{\end{env}}
    \newcommand{\emun}{\end{env}}
    \newcommand{\esp}{\end{env}}

\newcommand{\lf}{\vspace{2ex}}

\renewcommand{\bf}[1]{\textbf{#1}}
\renewcommand{\it}[1]{\textit{#1}}

\renewcommand{\sf}[1]{\textsf{#1}}

\renewcommand{\tt}[1]{\texttt{#1}}
\newcommand{\hl}[1]{\bf{\it{#1}}}

\newcommand{\mbf}[1]{\mathbf{#1}}

\newcommand{\cmc}[1]{\mathcal{#1}}
\newcommand{\eus}[1]{\mathscr{#1}}

\newcommand{\bb}[1]{\mathbb{#1}}
\newcommand{\nbd}[1]{$#1$\nobreakdash--}

\newcommand{\abs}[1]{\left\lvert#1\right\rvert}
\newcommand{\norm}[1]{\left\lVert#1\right\rVert}

\newcommand{\bfam}[1]{\bigl(#1\bigr)}
\newcommand{\Bfam}[1]{\Bigl(#1\Bigr)}
\newcommand{\AB}[1]{\langle#1\rangle}

\newcommand{\CB}[1]{\{#1\}}
\newcommand{\bCB}[1]{\bigl\{#1\bigr\}}
\newcommand{\BCB}[1]{\Bigl\{#1\Bigr\}}

\newcommand{\Matrix}[1]{\begin{pmatrix}#1\end{pmatrix}}

\newcommand{\set}[2][]{
    \ifthenelse{\equal{#1}{}}{
        \CB{#2}}{
        \CB{#1~|~#2}}}
\newcommand{\bset}[2][]{
    \ifthenelse{\equal{#1}{}}{
        \bCB{#2}}{
        \bCB{#1~|~#2}}}
\newcommand{\Bset}[2][]{
    \ifthenelse{\equal{#1}{}}{
        \BCB{#2}}{
        \BCB{#1~\big|~#2}}}

\DeclareMathOperator{\ls}{\normalfont\sf{span}}

\newcommand{\C}{\bb{C}}

\newcommand{\E}{\bb{E}}

\newcommand{\cA}{\cmc{A}}
\newcommand{\cB}{\cmc{B}}

\newcommand{\cM}{\cmc{M}}

\newcommand{\sB}{\eus{B}}

\newcommand{\U}{\mbf{1}}

\ifthenelse{\boolean{nosectionequn}}{}{
    \numberwithin{equation}{section}
    }

\ifthenelse{\boolean{article}\or\boolean{letter}\or\boolean{nosectionequn}}{
    \setboolean{nosectionappendix}{true}}{}
\ifthenelse{\boolean{nosectionappendix}}{}{
    \renewcommand{\appendix}{
        \chapter*{\appendixname}
        \addcontentsline{toc}{chapter}{\appendixname}
        \renewcommand{\thesection}{\Alph{section}}
        \setcounter{section}{0}}}
   
\ifthenelse{\boolean{report}\or\boolean{book}}{
    }{}

\ifthenelse{\boolean{notheorem}}{}{
        \newcommand{\definame}{Definition.}
        \newcommand{\propname}{Proposition.}
        \newcommand{\lemname}{Lemma.}
    \ifthenelse{\boolean{german}}{
        \newcommand{\exname}{Beispiel.}
        \newcommand{\remname}{Bemerkung.}
        \newcommand{\obname}{Beobachtung.}
        \newcommand{\thmname}{Satz.}
        \newcommand{\corname}{Korollar.}
        \newcommand{\proofname}{Beweis.}}{
        \newcommand{\exname}{Example.}
        \newcommand{\remname}{Remark.}
        \newcommand{\obname}{Observation.}
        \newcommand{\thmname}{Theorem.}
        \newcommand{\corname}{Corollary.}
        \newcommand{\proofname}{Proof.}}
    \theoremheaderfont{\normalfont\bfseries}
    \theoremstyle{change}
        \theorembodyfont{\rmfamily}
            \newtheorem{emp}{}[section]
                \newcommand{\bemp}[1][]{
                    \begin{emp}\hskip-\labelsep\bf{#1}\hskip\labelsep}
                \newcommand{\eemp}{\end{emp}}
\newtheorem{itemp}[emp]{}
                \newcommand{\bitemp}[1][]{
                    \begin{itemp}\hskip-\labelsep\bf{#1}\hskip\labelsep\normalfont\itshape}
                \newcommand{\eitemp}{\end{itemp}}
            \newtheorem{ex}[emp]{\exname}
                \newcommand{\bex}{\begin{ex}}
                \newcommand{\eex}{\end{ex}}
            \newtheorem{defi}[emp]{\definame}
                \newcommand{\bdefi}{\begin{defi}}
                \newcommand{\edefi}{\end{defi}}
            \newtheorem{rem}[emp]{\remname}
                \newcommand{\brem}{\begin{rem}}
                \newcommand{\erem}{\end{rem}}
            \newtheorem{ob}[emp]{\obname}
                \newcommand{\bob}{\begin{ob}}
                \newcommand{\eob}{\end{ob}}
        \theorembodyfont{\normalfont\itshape}
            \newtheorem{thm}[emp]{\thmname}
                \newcommand{\bthm}{\begin{thm}}
                \newcommand{\ethm}{\end{thm}}
            \newtheorem{prop}[emp]{\propname}
                \newcommand{\bprop}{\begin{prop}}
                \newcommand{\eprop}{\end{prop}}
            \newtheorem{cor}[emp]{\corname}
                \newcommand{\bcor}{\begin{cor}}
                \newcommand{\ecor}{\end{cor}}
            \newtheorem{lem}[emp]{\lemname}
                \newcommand{\blem}{\begin{lem}}
                \newcommand{\elem}{\end{lem}}
    \newcommand{\qedsymbol}{~\rule[-0.35mm]{2mm}{2mm}}
    \newcounter{proof}[emp]
    \newenvironment{Proof}[1]{
        \vspace{1ex}
        \renewcommand{\item}[1][\stepcounter{proof}(\roman{proof})]%
            {##1\hskip\labelsep}
        \noindent\textsc{#1\hskip\labelsep}}{
        \nolinebreak\qedsymbol}
    \newcommand{\proof}[1][\proofname]{
        \begin{Proof}{#1}\ignorespaces}
    \newcommand{\qed}{\end{Proof}}
    \newcommand{\noqed}{
        \renewcommand{\qedsymbol}{}
        \end{Proof}}}

    \usepackage{times,t1enc}
    \usepackage{mymathptm} 
    \usepackage[0.85]{fsscale}

\begin{document}

\title{Von Neumann Modules, Intertwiners and Self-Duality\thanks{This work is supported by DAAD and ISI Bangalore}}
\author{}
\author{
~\\
Michael Skeide\\\\
{\small\itshape Università degli Studi del Molise}\\
{\small\itshape Dipartimento S.E.G.e S.}\\
{\small\itshape Via de Sanctis}\\
{\small\itshape 86100 Campobasso, Italy}\\
{\small{\itshape E-mail: \tt{skeide@math.tu-cottbus.de}}}\\
{\small{\itshape Homepage: \tt{http://www.math.tu-cottbus.de/INSTITUT/lswas/\_skeide.html}}}\\
}
\date{Bangalore, July 2003}

{
\renewcommand{\baselinestretch}{1}
\maketitle

\vfill

\begin{abstract}
We apply the ideas of Muhly, Skeide and Solel \cite{MSS02p} of considering von Neumann $\cB$--mod\-ules as intertwiner spaces for representations of $\cB'$ to obtain a new, simple and self-contained proof for self-duality of von Neumann modules. This simplifies also the approach of \cite{MSS02p}.
\end{abstract}

}

\clearpage



\section{Introduction}

Let $E$ be a Hilbert module over a von Neumann algebra $\cB\subset\sB(G)$ acting (non-degenerately) on the Hilbert space $G$. We define the Hilbert space $H=E\odot G$ as the interior \hl{tensor product over $\cB$} of the right \nbd{\cB}module $E$ and the \nbd{\cB}\nbd{\C}module $G$ with inner product $\AB{x_1\odot g_1,x_2\odot g_2}=\AB{g_1,\AB{x_1,x_2}g_2}$. Every $x\in E$ gives rise to a mapping $L_x\colon g\mapsto x\odot g$ in $\sB(G,H)$ and it is easy to verify that $L{xb}=L_xb$ and $L_x^*L_y=\AB{x,y}$.

We, therefore, may and will identify every Hilbert \nbd{\cB}module over a von Neumann algebra $\cB\subset\sB(G)$ as a concrete submodule $E\subset\sB(G,H)$ of operators, where $H=E\odot G$. Following Skeide \cite{Ske00b,Ske01} we say $E$ is a \hl{von Neumann \nbd{\cB}module}, if $E$ is strongly closed in $\sB(G,H)$.

On $H$ we define a (normal unital) representation $\rho'\colon\cB'\rightarrow\sB(H)$ of the commutant $\cB'$ of $\cB$ by $\rho'(b')(id_E\odot b')$. (This is well-defined, because $b'$ is a bilinear mapping on the \nbd{\cB}\nbd{\C}module $G$, and also checking normality is routine.) In the special case when $E$ is the GNS-module of a completely positive mapping with values in $\cB$ (see Paschke \cite{Pas73}), $\rho'$ is known as \it{commutant lifting} (Arveson \cite{Arv69}).

Following Skeide \cite{Ske98,Ske00b}, the \hl{\nbd{\cB'}center} of the \nbd{\cB'}\nbd{\cB'}module $\sB(G,H)$ is defined as
\beqn{
C_{\cB'}(\sB(G,H))
~=~
\CB{x\in\sB(G,H)\colon\rho'(b')x=xb'~(b'\in\cB')}.
}\eeqn
As observed, for instance, by Goswami and Sinha \cite{GoSi99}, it is easy to check that $C_{\cB'}(\sB(G,H))$ is a itself a von Neumann \nbd{\cB}module.

Clearly, $E$ is contained in $C_{\cB'}(\sB(G,H))$.

It is the starting point in \cite{MSS02p} to show that $E$ is all of $C_{\cB'}(\sB(G,H))$. Once known that von Neumann modules are \hl{self dual }, i.e.\ every bounded right linear mapping $\Phi\colon E\rightarrow\cB$ (so-called \hl{\nbd{\cB}functionals}) has the form $\AB{x,\bullet}$ for a (unique) $x\in E$, (see \cite{Ske00b,Ske01} for proof using \it{complete quasi orthonormal systems}, a suitable geralization of orthonormal bases in Hilbert spaces) this task is easy: Like for Hilbert spaces a strongly closed (and, therefore, self-dual) submodule with zero-complement is all. And since $EG$ is total in $H$ the complement of $E$ in $C_{\cB'}(\sB(G,H))$ is, indeed, $\CB{0}$.

\brem
The other important observation in \cite{MSS02p} is that for an arbitrary (normal unital) representation $\rho'$ of $\cB'$ on a Hilbert space, $C_{\cB'}(\sB(G,H))$ is a von Neumann \nbd{\cB}module acting totally on $G$, what gives a one-to-one correspondence between von Neumann \nbd{\cB}modules contained in $\sB(G,H)$ and representations $\rho'$ of $\cB'$ on $H$.

This approach becomes particularly fruitful, when the von Neumann modules are two-sided so that there is arround another (normal unital) representation $\rho$ on $H$ of a second von Neumann algebra $\cA$. Switching the roles of $\cB$ and $\cA'$, the result is a duality between \nbd{\cA}\nbd{\cB}modules and \nbd{\cB'}\nbd{\cA'}modules generalizing the duality between a von Neumann algebra and its commutant. One application is a complete theory of normal representations of the adjointable operators on a von Neumann \nbd{\cB}module on a von Neumann \nbd{\cA}module (this can be, e.g., a Hilbert space).
\erem

In this short note we show that $E=C_{\cB'}(\sB(G,H))$ and we show that $C_{\cB'}(\sB(G,H))$ is self-dual, thus, showing that $E$ is self-dual. The only prerequisit for the first statement is von Neumann's \it{double commutant theorem}, the only prerequisit for the second statement is a technical lemma which asserts that every \nbd{\cB}functional $\Phi$ is represented by an operator in $\sB(H,G)$ (see below). The simplicity of the proofs improves also accessibility of \cite{MSS02p} and, therefore, of the whole theory of von Neumann modules.

\section{$E=C_{\cB'}(\sB(G,H))$}

Every $a\in\sB^a(E)$ (the \hl{algebra of adjointable operators} on $E$) gives rise to a bounded operator $x\odot g\mapsto ax\odot g$ on $H$. In that way, we identify $\sB^a(E)$ as a \nbd{*}subalgebra of $\sB(H)$. It is easy to see that $\sB^a(E)$ is a von Neumann subalgebra of $\sB(H)$.

It follows that the \it{matrix \nbd{*}algebra}
\beqn{
\cM
~=~
\Matrix{\cB&E^*\\E&\sB^a(E)}
}\eeqn
with the obvious operations is a von Neumann algebra on $G\oplus H$. Let us compute its commutant.

\bprop
The commutant of $\cM$ is $\cM'=\BCB{\Bfam{\substack{b'\\\\0}~\substack{0\\\\\rho'(b')}}\colon b'\in\cB'}$.
\eprop

\proof
Suppose $\Bfam{\substack{b'\\\\x'}~\substack{{y'}^*\\\\a'}}\in\sB(G\oplus H)$ is an element in $\cM'$. As it must commute with $\Bfam{\substack{b\\\\0}~\substack{0\\\\0}}$ $(b\in\cB)$ we find
\beqn{
\Matrix{b'b&0\\x'b&0}
~=~
\Matrix{bb'&b{y'}^*\\0&0}.
}\eeqn
As this must hold for all $b\in\cB$ (in particular also for $b=\U$), we find $x'=y'=0$ and $b'\in\cB'$. The remaining part $\Bfam{\substack{b'\\\\0}~\substack{0\\\\a'}}$ must commute with $\Bfam{\substack{0\\\\x}~\substack{0\\\\0}}$ $(x\in E)$. Therefore,
\beqn{
\Matrix{0&a'x\\0&0}
~=~
\Matrix{0&xb'\\0&0}.
}\eeqn
We find $a'(x\odot g)=a'xg=xb'g=\rho'(b')(x\odot g)$ for all $x\in E$, $g\in G$ and, therefore, $a'=\rho'(b')$.\qed

\lf
The commutant of $\cM'$ is, clearly,
\beqn{
\cM''
~=~
\Matrix{\cB&C_{\cB'}(\sB(H,G))\\C_{\cB'}(\sB(G,H))&\rho'(\cB')'}.
}\eeqn
By the \it{double commutant theorem} $\cM''=\cM$. Therefore, we do not only show the statement of this section's headline, but, as an additional benefit, we identify also $\sB^a(E)$ as the commutant of the image of $\cB'$ under $\rho'$. (This can also be done by using \it{Morita equivalence} for von Neumann algebras; see Rieffel \cite{Rie74a}.)

\bprop
$E=C_{\cB'}(\sB(G,H))$ and $\sB^a(E)=\rho'(\cB')'$.
\eprop

\section{$C_{\cB'}(\sB(G,H))$ is self-dual}

A \nbd{\cB}functional $\Phi\in\cB^r\bfam{C_{\cB'}(\sB(G,H)),\cB}$ gives rise to a linear mapping
\beqn{
L_\Phi\colon
\ls C_{\cB'}(\sB(G,H))G
~\rightarrow~
G
~~~~~~~~~~~
L_\Phi(x\odot g)
~=~
(\Phi x)g.
}\eeqn
The proof of the following lemma consists, essentially, in showing that for computing the operator norm of $L_\Phi$ it is sufficient to take the supremum only over elementary tensors $x\odot g$ ($\norm{x}\le1,\norm{g}\le1$).

\blem\label{boundlem}
$\norm{L_\Phi}=\norm{\Phi}$. Therefore, $L_\Phi$ extends to a bounded operator in $\sB(H,G)$ identified with $\Phi$.
\elem

\proof
(Sketch only. See \cite{Ske00b,Ske01} for details.) Suppose that there is a cyclic vector $g_0\in G$, i.e.\ $\cB g_0$ is dense in $G$. (Otherwise, use a decomposition of $G$ into subspaces $G_\alpha$ cyclic for $\cB$ and take into account the facts, firstly, that also $H$ decomposes accordingly into cyclic subspaces $H_\alpha$ and, secondly, that the norm of an element in a direct sum of operator spaces $\sB(G_\alpha,H_\alpha)$ is just the supremum over the single norms.) Then every element in $H$ can be approximated by those of the form $h=x\odot g_0$. Use polar decomposition $x=x_0\abs{x}$ of $x$ and put $g=\abs{x}g_0$. Then $\norm{h}=\norm{g}$ because  $g\in\abs{x}G$. In particular, every unit vector in $H$ can be approximated by $x\odot g$ where $x$ is a partial isometry in $E$ and $g$ is a unit vector in $G$.\qed

\bprop\label{sdprop}
$C_{\cB'}(\sB(G,H))$ is self-dual.
\eprop

\proof
From $\Phi\rho'(b')(x\odot g)=\Phi(x\odot b'g)=\Phi xb'g=b'\Phi xg=b'\Phi(x\odot g)$ we see that $\Phi$ intertwines $b'$ and $\rho'(b')$. Therefore the adjoint $y=\Phi^*$ of $\Phi$ is an element in $C_{\cB'}(\sB(G,H))$ such that $\Phi x=\AB{y,x}$ for all $x\in C_{\cB'}(\sB(G,H))$.\qed

\section{Synthesis}

\bthm
Every von Neumann \nbd{\cB}module is self-dual.
\ethm

\proof
$E=C_{\cB'}(\sB(G,H))$ and $C_{\cB'}(\sB(G,H))$ is self-dual.\qed

\brem
It seems that Lemma \ref{boundlem} forms always an essential part of the proofs of self-dualtiy, which cannot be replaced by simpler arguments.
\erem

\brem
To be honest, we should mention that (unlike the approach by complete quasi orthonormal systems) the preceding arguments cannot be used to show the \it{Riesz representation theorem} (Hilbert spaces are self-dual), but, actually, reduce the statement about von Neumann modules to that about Hilbert spaces. An equivalent form of the \it{Riesz representation theorem} is that all bounded operators between Hilbert spaces have an adjoint. Without this, in the proof of Proposition \ref{sdprop} it was not possible to pass from $\Phi$ to $\Phi^*$. One may see the failure of the argument clearly, by taking $\cB=\C$ and $G=\C$ and for $H$ only a pre-Hilbert space. This is, actually, the only place in these notes, where we are not able to write down an adjoint explicitly on the algebraic domain $\ls EG$. (The adjoint of $x\colon g\mapsto x\odot g$ is, of course, $x^*\colon y\odot g\mapsto\AB{x,y}g$.)
\erem

\baselineskip2.8ex plus 0ex minus 0.1ex 

\lf\noindent
\bf{Acknowledgement.}~
These notes were written during a two months stay of the author at ISI Bangalore. The author wishes to express his gratitude to Prof.\ B.V.R.\ Bhat and the ISI for warm hospitality and to the DAAD for travel support.


\newcommand{\Swap}[2]{#2#1}\newcommand{\Sort}[1]{}
\providecommand{\bysame}{\leavevmode\hbox to3em{\hrulefill}\thinspace}
\providecommand{\MR}{\relax\ifhmode\unskip\space\fi MR }
\providecommand{\MRhref}[2]{%
  \href{http://www.ams.org/mathscinet-getitem?mr=#1}{#2}
}
\providecommand{\href}[2]{#2}


\end{document}